\documentclass[12pt]{article}

\usepackage{sola}

\usepackage{pgf,ytableau}

\usepackage{algorithm}
\usepackage[noend]{algpseudocode}
\usepackage{xcolor}

\def\addtab#1={#1\;&=}

\def\meeq#1{\def\ccr{\\\addtab}
 \begin{align*}
 \addtab#1
 \end{align*}
  }  
  
  \def\leqaddtab#1\leq{#1\;&\leq}

\def\vc#1{\mbox{\boldmath$#1$\unboldmath}}

\def\pr(#1){\left({#1}\right)}
\def\br[#1]{\left[{#1}\right]}
\def\fbr[#1]{\!\left[{#1}\right]}
\def\set#1{\left\{{#1}\right\}}
\def\ip<#1>{\left\langle{#1}\right\rangle}
\def\iip<#1>{\left\langle\!\langle{#1}\right\rangle\!\rangle}

\def\fpr(#1){\!\pr({#1})}

\def\mapengine#1,#2.{\mapfunction{#1}\ifx\void#2\else\mapengine #2.\fi }

\def\map[#1]{\mapengine #1,\void.}

\def\mapenginesep_#1#2,#3.{\mapfunction{#2}\ifx\void#3\else#1\mapengine #3.\fi }

\def\mapsep_#1[#2]{\mapenginesep_{#1}#2,\void.}

\def\vcbr{\br}

\def\bvect[#1,#2]{
{
\def\dots{\cdots}
\def\mapfunction##1{\ | \  ##1}
	\sopmatrix{
		 \,#1\map[#2]\,
	}
}
}

\def\vect[#1]{
{\def\dots{\ldots}
	\vcbr[{#1}]
}}

\def\vectt[#1]{
{\def\dots{\ldots}
	\vect[{#1}]^{\top}
}}

\def\Vectt[#1]{
{
\def\mapfunction##1{##1 \cr} 
\def\dots{\vdots}
	\begin{pmatrix}
		\map[#1]
	\end{pmatrix}
}}

\def\tF_#1{{\tt F}_{#1}}

\def\tFC_#1{{\tt T}_{#1}}

\def\Problem#1#2\par{\begin{problem}\label{Problem:#1} \prooflabel{#1}#2\end{problem}}
\def\Theorem#1#2\par{\begin{theorem}\label{Theorem:#1} \prooflabel{#1}#2\end{theorem}}
\def\Conjecture#1#2\par{\begin{conjecture}\label{Conjecture:#1} \prooflabel{#1}#2\end{conjecture}}
\def\Proposition#1#2\par{\begin{proposition}\label{Proposition:#1} \prooflabel{#1}#2\end{proposition}}
\def\Definition#1#2\par{\begin{definition}\label{Definition:#1} \prooflabel{#1}#2\end{definition}}
\def\Corollary#1#2\par{\begin{corollary}\label{Corollary:#1} \prooflabel{#1}#2\end{corollary}}
\def\Lemma#1#2\par{\begin{lemma}\label{Lemma:#1} \prooflabel{#1}#2\end{lemma}}
\def\Example#1#2\par{\begin{example}\label{Exercise:#1} \prooflabel{#1}#2\end{example}}
\def\Remark #1\par{\begin{remark*} #1\end{remark*}}

\def\Section#1#2.{\section{#2}\label{Section:#1} \prooflabel{#1}}

\def\secref#1{Section~\ref{Section:#1}}

\def\lmref#1{Lemma~\ref{Lemma:#1}}
\def\propref#1{Proposition~\ref{Proposition:#1}}

\def\thref#1{Theorem~\ref{Theorem:#1}}
\def\defref#1{Definition~\ref{Definition:#1}}
\def\probref#1{Problem~\ref{Problem:#1}}
\def\corref#1{Corollary~\ref{Corollary:#1}}

\def\Ai{{\rm Ai}\,}

\def\elllRpz_#1{\ell_{#1{\rm z}}^{(\lambda,R),p}}

\def\sopmatrix#1{\begin{pmatrix}#1\end{pmatrix}}

\newtheorem{problem}{Problem}
\newtheorem{lemma}{Lemma}
\newenvironment{remark*}{{\it Remark} }{}
\newtheorem{theorem}{Theorem}
\newtheorem{conjecture}{Conjecture}
\newtheorem{corollary}{Corollary}
\newtheorem{proposition}{Proposition}

\theoremstyle{definition}
\newtheorem{definition}{Definition}
\newtheorem{example}{Example}

\begin{document}




\refmodetrue

\def\rhos{\rho^{\rm s}}
\def\rhoe{\rho^{\rm e}}
\def\GZn{\hbox{GZ}(n)}
\def\Spec{\hbox{Spec}}
\def\Cont{\hbox{\rm Cont}}
\def\cont{\hbox{\rm cont}}
\def\kron{\hbox{\rm kron}}
\def\vec{\hbox{\rm vec}\,}
\def\mspan{\hbox{\rm span}\,}
\def\colspan{\hbox{\rm colspan}\,}
\def\GL{\hbox{\rm GL}}
\def\mdiag{\hbox{\rm diag}}
\def\Y{{\mathcal Y}}
\def\algref#1{Algorithm~\ref{Algorithm:#1}}
\def\bfc{{\bf c}}
\def\bfx{{\bf x}}

\def\addition{\relax}

\def\probsref[#1,#2]{Problems \ref{Problem:#1} and \ref{Problem:#2}}

\def\bbR{{\mathbb R}}
\def\bbZ{{\mathbb Z}}
\def\bbC{{\mathbb C}}
\def\bbP{{\mathbb P}}
\def\bfP{{\mathbf P}}
\def\bfQ{{\mathbf Q}}
\def\bfe{\vc e}
\def\bfq{\vc q}
\def\bfu{\vc u}

\def\Orth{\hbox{\rm O}}
\def\Ord{{\mathcal O}}

\def\hexnumber#1{\ifcase#1 0\or1\or2\or3\or4\or5\or6\or7\or8\or9\or
 A\or B\or C\or D\or E\or F\fi}
\edef\msbhx{\hexnumber\symAMSb}   
\mathchardef\emptyset="0\msbhx3F

\authord={Sheehan Olver\footnote{Department of Mathematics,  Imperial College London, UK, {\tt s.olver@imperial.ac.uk}}}

\titled={Representations of the symmetric group are  decomposable in polynomial time}

\maketitle

\Abstract
We introduce an algorithm to decompose orthogonal matrix representations of the symmetric group over the reals into irreducible representations, which as a by-product also computes the multiplicities of the irreducible representations. The algorithm applied to a $d$-dimensional representation of $S_n$  is shown to have a complexity of $\Ord(n^2 d^3)$ operations for determining {\addition which irreducible representations are present and their corresponding} multiplicities and a further $\Ord(n d^4)$ operations to fully decompose representations with non-trivial multiplicities.  These complexity bounds are pessimistic and in a practical implementation using floating point arithmetic and exploiting sparsity we observe better complexity. We demonstrate this algorithm on the problem of computing multiplicities of two tensor products of irreducible representations (the Kronecker coefficients problem) as well as higher order tensor products. For hook and hook-like irreducible representations the algorithm  has polynomial complexity as $n$ increases. {\addition We also demonstrate an application to constructing a basis of multivariate orthogonal polynomials with respect to a tensor product weight so that applying a permutation of variables induces an irreducible representation.}

	
\Section{intro}  Introduction.


Let $\rho : S_n \rightarrow \Orth(d)$ be an {\it orthogonal matrix representation} of the symmetric group, acting on the vector space $V = \bbR^d$.  In other words $\rho$ is a {\it homomorphism}: for all $a,b \in S_n$ we have $\rho(ab) = \rho(a)\rho(b)$. We consider the problem of decomposing the representation into {\it irreducible representations}, that is decomposing
$$
	V \cong V_1^{\oplus a_1} \oplus \cdots \oplus V_r^{\oplus a_r}
$$
where $\rho_1, \ldots, \rho_r$ are each irreducible representations acting on the vector spaces $V_1, \ldots, V_r$, {\addition where in the setting of the symmetric group these are each associated with a different partition of $n$}. 
 We  assume that we are given the generators of a representation of $S_n$ as symmetric (which must also be orthogonal) matrices $\rho(\tau_1), \ldots, \rho(\tau_{n-1}) \in \Orth(d)$ where $\tau_\ell = (\ell\ (\ell+1))$ are the simple transpositions, i.e., the Coxeter generators. 

 {\addition
Given generators $\rho(\tau_1), \ldots, \rho(\tau_{n-1}) \in \Orth(d)$ of  an orthogonal matrix representation $\rho :  S_n \rightarrow \Orth(d)$, we shall first consider the following fundamental problems:

 \Problem{present}  Determine the irreducible representations present in the above decomposition.
 
 \Problem{multiplicities}   Determine the (non-zero) {\it multiplicities} $a_1,\ldots,a_r$ of each irreducible representation, that is, how many times does each irreducible representation occur.
 
 We will then  solve the question of decomposing a representation by solving an equivalent linear algebra problem:
 

\Problem{blockdiag}
Compute an orthogonal matrix $Q \in \Orth(d)$ that block-diagonalises the representation  so that
$$
Q^\top \rho(g) Q = \sopmatrix{ \rho_1(g)^{\oplus a_1} \\ & \ddots \\ && \rho_r(g)^{\oplus a_r}}
$$
where $\rho_1,\ldots,\rho_r$ are distinct  irreducible representations and $a_1,\ldots,a_r$ are their multiplicities. Here we use the notation
$$
\rho_k(g)^{\oplus a_k} = I_{a_k} \otimes \rho_k(g) =  \sopmatrix{ \rho_k(g) \\ & \ddots \\ && \rho_k(g)} \in \Orth(d_k a_k)
$$
where $\rho_k(g) \in \Orth(d_k)$ and $d_k$ is the dimension of the irreducible representation.

In this paper we will prove that both problems are computable (assuming real arithmetic). In particular \probsref[present,multiplicities] can be solved via \algref{multiplicities} in $\Ord(n^2 d^3)$ operations (\corref{prob1})  and \probref{blockdiag} can be solved via \algref{blockdiag} in a further $\Ord(n d^4)$ operations (\corref{fullyreduce}).   Going further, the practical realisation of \algref{blockdiag} appears to achieve $\Ord(n d^3)$ operations by exploiting sparsity, but this more efficient complexity result is not yet proven, see discussion in \secref{blockdiag}.

\bigskip

Prior work on these problems begins with Dixon \cite{DixonComputeReps}, who constructed an iterative algorithm for reducing representations of general finite groups to irreducible components, which can be used to solve  \probsref[present,multiplicities]  in the limit.  However, it does not consider the decomposition of repeated copies of irreducible representations hence does not solve \probref{blockdiag}, that is, it only computes the {\it canonical representation} as defined in \cite{serre1977linear}. Furthermore, the algorithm requires the computation of eigenvectors and eigenvalues which also depends on iterative algorithms which may fail.    There is also related work on decomposing representations of compact groups \cite{rosset2021replab}.

It is possible to solve \probref{multiplicities} by using the orthogonality relationships of the characters $\chi(g) := {\rm Tr} \rho(g)$ with respect to the inner product
$$
\langle \phi, \psi \rangle_{S_n} := {1 \over n!} \sum_{g \in S_n} \phi(g) \psi(g),
$$
see for example \cite[Theorem 4]{serre1977linear}. As characters are class functions they depend only on the conjugacy class of the permutation, therefore this sum can be reduced to one over all partitions, whose number are given by the partition function $p(n)$. The complexity of the resulting algorithm depends on the cost in creating a representation $\rho(g)$ for $g$ running through the a set of permutations in the given conjugacy class. If these can be constructed directly in $\Ord(d)$ operations then the total complexity is $\Ord(p(n) d)$. However, if one is only given generators (that may be dense matrices)  we need to construct other representations by multiplying as many as $n$ matrices, hence the total complexity is $\Ord(p(n) n d^{2 + \kappa})$ where the $\kappa$ depends on which matrix-multiplication algorithm is used ($\kappa = 1$ with the standard definition of matrix multiplication but can be improved using fast matrix multiplication). In either case, as $p(n)$ grows faster than algebraically it does not result in a polynomial complexity algorithm. 

In other work, Serre \cite[Theorem 8, Proposition 8]{serre1977linear} introduced explicit projections to the invariant vector spaces corresponding to irreducible representations, and the column rank of these projectors can be used to solve \probref{multiplicities}. However, the projections require computing sums over every element of the group, that is,  in the case of the symmetric group the complexity for constructing each projector is $\Ord(n! d^2)$ operations. The number of projectors needed is $\Ord(d^2)$ so the total complexity is $\Ord(n! d^4)$ operations.  Furthermore, it does not give a fast approach to solve \probref{present}, that is,  it cannot determine which irreducible representations are present (i.e., have non-zero multiplicity) without potentially checking all $p(n)$ possibilities which grows super-algebraically with $n$. Finally, whilst a QR decomposition of the projectors can be used to construct an orthogonal matrix $Q$ that block-diagonalises a representation into irreducible representation there is no guarantee that the irreducible representations in the same isomorphism class will be identical, that is, one likely needs to apply \algref{blockdiag} (or an equivalent procedure) to enforce these representations are identical.

The complexity of Serre's approach for constructing projectors can be dramatically reduced to polynomial in $n$ using the ideas introduced by Hymabaccus and Pasechnik in  \cite{hymabaccus2020decomposing}, see also the related software \cite{hymabaccus2020repndecomp}. In particular, in the case of the symmetric group the complexity of computing the projectors is reducible to $\Ord(n^2 d^2)$ operations for a total of $\Ord(n^2 d^4)$ operations to solve \probref{multiplicities}, assuming one has computed which irreducible represents are present (i.e. solved \probref{present}). Unfortunately this work does not provide a fast way to solve \probref{present} and hence the total cost is still not polynomial, unless combined with the proposed \algref{multiplicities}. 

\Remark The work of \cite{hymabaccus2020decomposing,hymabaccus2020repndecomp} takes a more computational algebra approach and has the significant benefit that the computations can be performed exactly. As the algorithms in this paper work with reals and compute orthogonal matrices that cannot in general be represented exactly on a computer  the practical implementation inevitably requires floating point arithmetic operations, hence the calculations will not be exact due to round-off errors. However, as the algorithms are built on top of  well-understand numerical linear algebra algorithms such as those for computing eigenvalues of symmetric matrices the computations are reliable, as demonstrated in \figref{cond_p}. Going a step further techniques from Validated Numerics may yield a rigorous implementation or verification of the results, a topic we discuss in \secref{future}.

}

Solving \probref{multiplicities} in the special case where the representation is a tensor product of irreducible representations of $S_n$ is one approach to computing Kronecker coefficients, a problem that is known to be NP-hard  (in particular, $\#$P-hard \cite{ikenmeyer2017vanishing,burgisser2008complexity}). In certain settings where the corresponding partitions have a fixed length there is a combinatorial algorithm that can compute specific Kronecker coefficients in polynomial time  \cite{christandl2012computing,pak2017complexity}. It is difficult to make a direct comparison between our complexity results when applied to the Kronecker coefficients problem:  our complexity results depend on the dimension of the corresponding irreducible representations and we obtain all non-zero multiplicities without having to deduce the zero multiplicities, whereas   the complexity results of \cite{pak2017complexity} largely depend on the length of the corresponding partitions (including that of the partition corresponding to the irreducible representation whose multiplicity is being computed) and give no approach to deduce which multiplicities are non-zero apart from testing all possible $p(n)$ partitions, whose number grows faster than polynomial in $n$.

\bigskip

{\addition We will now construct a concrete example that will be used throughout this paper. Denote the representation coming from permutation matrices as  $\rhos_n : S_n \rightarrow  \bbR^{n \times n}$, where $\rhos_n(g)$ is the identity matrix $I_n$ with the rows permuted according to $g$.   It has symmetric generators $\rhos_n(\tau_\ell)$ with the $\ell$ and $(\ell+1)$th rows permuted. We will consider the representation  $\rhoe : S_4 \rightarrow O(6)$ defined by:
\begin{align*}
\rhoe(\tau_\ell) &:= \begin{pmatrix} \rhos_4(\tau_\ell) \\ & \rhos_2(\tau_1) \end{pmatrix} \qquad \hbox{e.g.} \qquad \rhoe(\tau_2) = \begin{pmatrix} 1 \\ & 0 & 1 \\ & 1 & 0 \\ &&& 1 \\ &&&& 0 & 1 \\ &&&& 1 & 0
\end{pmatrix}
\end{align*}
This representation is block-diagonalised  by the matrix
\begin{equation}\label{eq:permutationQ}
Q = \sopmatrix{
		0 	& 1/\sqrt{2} 	& 1/\sqrt{6} 	& 	1/\sqrt{12}		&	1/2	&	0 \cr
                 0	& -1/\sqrt{2} 	& 1/\sqrt{6}	&	1/\sqrt{12} 	& 1/2 	& 	0 \cr
                 0	&0		 	& -\sqrt{2/3}	&	1/\sqrt{12} 	& 1/2 	& 	0 \cr
                 0	&0			& 0			&	-\sqrt{3}/2 		& 1/2		&	0 \\ 
                 1/\sqrt{2} &0		& 0			&	0			&0		& 1\/\sqrt{2} \\
                 -1/\sqrt{2} &0		& 0			&	0			& 0		& 1/\sqrt{2} 
                 }
\end{equation}
}
A quick check shows that we have block-diagonalised the generators (and hence the representation) into {\addition a} $3 \times 3$ and {\addition three} $1 \times 1$ sub-blocks:{\addition
\meeq{
Q^\top \rhoe(\tau_1) Q =
\sopmatrix{-1 &\vline & &   & & \vline & &\vline &   \\ 
\hline
& \vline & -1  &  &  & \vline & &\vline &\\ 
& \vline &  & 1 &   & \vline &  &\vline &\\ 
& \vline &  &  & 1 &\vline & & \vline &\\ 
\hline
& \vline &  &  &   & \vline &  1&\vline &\\ 
\hline
& \vline &  &  &   & \vline &  &\vline &1 \\ 
}, \ccr
Q^\top \rhoe(\tau_2) Q = \sopmatrix{-1 &\vline & &   & & \vline & &\vline &   \\ 
\hline
& \vline & 1/2  & \sqrt{3}/2 &  & \vline & &\vline &\\ 
& \vline & \sqrt{3}/2  & -1/2 &   & \vline &  &\vline &\\ 
& \vline &  &  & 1 &\vline & & \vline &\\ 
\hline
& \vline &  &  &   & \vline &  1&\vline &\\ 
\hline
& \vline &  &  &   & \vline &  &\vline &1 \\ 
}, \ccr
Q^\top\rhoe(\tau_3) Q =
\sopmatrix{-1 &\vline & &   & & \vline & &\vline &   \\ 
\hline
& \vline & 1  &  &  & \vline & &\vline &\\ 
& \vline &  & 1/3 &  2\sqrt{2}/3  & \vline &  &\vline &\\ 
& \vline &  & 2\sqrt{2}/3  &- 1/3 &\vline & & \vline &\\ 
\hline
& \vline &  &  &   & \vline &  1&\vline &\\ 
\hline
& \vline &  &  &   & \vline &  &\vline &1 \\ 
}
}
}
Note that the sub-blocks are necessarily also representations of $S_4$.  That is, we have decomposed the permutation matrix representation
into {\addition three irreducible representations: one associated with the partition $4 = 3 +1$,  one associated with the trivial partition $4 = 4$ (with multiplicity two) and one associated with the sign partition $4 = 1 + 1 + 1 + 1$.} The algorithms we introduce computes these partitions and multiplicities as well as the matrix $Q$ from the generators $\rhoe(\tau_1),  \rhoe(\tau_2), \rhoe(\tau_3)$.

\bigskip


%
%
%
%
%
%
%
The paper is structured as follows:

{\noindent \secref{irreps}}: We review basics of representation theory of the symmetric group, including irreducible representations, their generators, the Gelfand--Tsetlin (GZ) algebra and its spectral properties for irreducible representations.

{\noindent \secref{rep2la}}: We detail how the problem of reducing an orthogonal matrix representation can be recast {\addition in terms of linear algebra, in particular} as \probref{blockdiag}.

{\noindent \secref{mult}}: We discuss how the GZ algebra can be used to solve \probref{multiplicities} via a joint spectrum problem involving commuting symmetric matrices. This also gives a way to reduce a representation into irreducible representations but cannot distinguish multiple copies of the same irreducible representation, that is, we {\addition only compute reduction to a} canonical representation.

{\noindent \secref{blockdiag}}: We discuss how irreducible representations involving multiple copies of the same irreducible representation can be {\addition fully} reduced, leading to the solution of \probref{blockdiag}.

{\noindent  \secref{algorithm}}: We outline the algorithms for solving \probref{blockdiag} and \probref{multiplicities} and discuss briefly their practical implementation using floating point arithmetic.

{\noindent  \secref{examples}}: We show examples including tensor products of irreducible representations (the Kronecker coefficients) and higher order analogues, using a floating-point arithmetic implementation of the proposed algorithm. We also demonstrate polynomial complexity for hook and hook-like irreducible representations. {\addition Finally, we consider the problem of constructing a basis of orthogonal polynomials with respect to a tensor product weight in $n$-dimensions so that permuting variables give rise to irreducible representations of $S_n$. 

{\noindent \secref{future}}: We briefly discuss potential future work including better complexity algorithms, adaptation to Coxeter groups, and applications to sparse discretisations arising in numerical quadrature and solutions to partial differential equations on geometries with symmetries.

\bigskip

 \Remark Most of the theoretical results can be adapted to non-orthogonal representations $\rho : S_n \rightarrow {\rm GL}_n(\bbR)$ or $\rho : S_n \rightarrow {\rm GL}_n(\bbC)$ though possibly with worse complexities. For example, the tridiagonalisation procedure in \lmref{cubicjointspec} only applies for symmetric matrices and hence computing the relevant nullspaces will have worse complexity. Moreover, the practical implementation with non-orthogonal representations may be less-reliable as algorithms involving non-orthogonal and non-symmetric matrices are prone to issues due to ill-conditioning, or in the case of computing an eigendecomposition may potentially fail.

}

\bigskip
\noindent{\bf Acknowledgments}: I thank Oded Yacobi (U. Sydney) for significant help in understanding the basics of representation theory and in particular \cite{Okounkov}, as well as Peter Olver (U. Minnesota), Alex Townsend (Cornell), and Marcus Webb (U. Manchester) for helpful suggestions on drafts.  {\addition We also thank the anonymous referees for their very helpful feedback.} This work was completed with the support of the EPSRC grant EP/T022132/1 “Spectral element methods for fractional differential equations, with applications in applied analysis and medical imaging” and the Leverhulme Trust Research Project Grant RPG-2019-144 “Constructive approximation
theory on and inside algebraic curves and surfaces”.

\Section{irreps} Irreducible representations of the symmetric group.

In this section we review some basic facts of representation theory of the symmetric group, roughly following \cite{Okounkov}. The irreducible representations can be identified with {\it partitions}:

\Definition{Partitions} A {\it partition of $n$},  is a tuple $\lambda = (\lambda_1,\ldots,\lambda_k)$ of integers $\lambda_1 \geq \cdots \geq \lambda_k \geq 1$ such that $n = \lambda_1 + \cdots + \lambda_k$. We use the notation $\lambda  \vdash n$ to denote that $\lambda$ is a partition of $n$.

A basis for the vector space associated with an irreducible representation can be identified with {\it Young tableaux}:

\Definition{YoungTableau} A {\it Young tableau} is  a chain of partitions $(\lambda^1,\ldots,\lambda^n)$ where $\lambda^j \vdash j$ such that $\lambda^{j+1}$ is equivalent to $\lambda^j$ with
one entry increased or one additional entry {\addition equal to one}.  The set of all Young tableaux of length $n$ is denoted $\Y_n$.

 A natural way to visualise a partition is via a {\it Young diagram}: if $\lambda = (\lambda_1,\ldots,\lambda_k)$ then a Young diagram consists of $k$ rows of boxes where the $j$th row has exactly $\lambda_j$ boxes. A natural way to visualise a Young tableau is by filling in the boxes of a Young Diagram  according to the
 order in which the new boxes appear in the sequence of Young Diagrams corresponding to $\lambda^j$.  For example, the Young tableau $((1), (2), (2,1), (2,1,1), (2,2,1), (3,2,1), (3,2,1,1))$ can be depicted
\begin{equation}\label{eq:ytableauexample}
 \begin{ytableau}
 1 & 2 & 6 \\
 3 & 5 \\
 4 \\
 7
 \end{ytableau}
\end{equation}

The number of Young tableaux corresponding to a given partition $\lambda \vdash n$ is called the {\it hook length} $d_\lambda$. Young tableaux can be used
to build an explicit irreducible matrix representation:

\Definition{irreps} Associated with any partition $\lambda$ is a {\it canonical orthogonal matrix irreducible representation} $\rho_\lambda : S_n \rightarrow \Orth(d_\lambda)$. The dimension of the irreducible representation
is the number of Young tableaux $d_\lambda$ and hence we can parameterise $\bbR^{d_\lambda}$ by an enumeration of Young tableaux $\set{y_1, \ldots, y_{d_\lambda}} = \Y_n$. The entries of the (symmetric) generators $\rho_\lambda(\tau_\ell)$ are
given by the following rules (cf.~\cite[(6.5)]{Okounkov}):
\begin{enumerate}
    \item   For the Young tableau $y_j$, if box $\ell$ and $\ell+1$ are in the same row then $\vc e_j^\top \rho_\lambda(\tau_\ell) \vc e_j = 1$.
    \item   For the Young tableau $y_j$, if box $\ell$ and $\ell+1$ are in the same column then $\vc e_j^\top \rho_\lambda(\tau_\ell) \vc e_j = -1$.    
    \item   For the Young tableau $y_j$, if box $\ell$ and $\ell+1$ are not in the same row or column then let $k$ be such that $y_k$
    is equal to $y_j$ except with boxes $\ell$ and $\ell+1$ swapped. If $(i_\ell,j_\ell)$ is the index of box $\ell$, then for the {\it axial distance} $r = j_{\ell+1}+i_\ell-i_{\ell+1} -j_\ell$ we have
    \meeq{
   \vc e_j^\top  \rho_\lambda(\tau_\ell) \vc e_j = 1/r,  \ccr
   \vc e_k^\top \rho_\lambda(\tau_\ell) \vc e_j  = \sqrt{1-1/r^2}
}
    \item All other entries are zero.
\end{enumerate}

%

We now introduce the Gelfand--Tsetlin (GZ) algebra, which is a commutative algebra that we shall utilise to determine the multiplicities of the irreducible representations as well as decompose a representation into its canonical representation.

\Definition{YJMgen} The {\it Young--Jucys--Murphy (YJM)-generators}  are
\meeq{
X_1 = 0, X_2 = (1\ 2), X_3 = (1\ 3) +  (2\ 3),\ldots,X_n = (1\ n) + \cdots + (n-1\ n)
}
where $(i\ j) \in S_n$ is cyclic notation for the permutation that swaps $i$ and $j$.

The Gelfand--Tsetlin (GZ) algebra, a sub-algebra of the group algebra of $S_n$,  is generated by the YJM generators \cite[Corollary 2.6]{Okounkov}:
$
\GZn = \ip< X_1, \ldots, X_n >.
$
Note that a representation  $\rho : S_n \rightarrow  \Orth(d)$ induces a representation $ \rho : \GZn \rightarrow \bbR^{d \times d}$. We shall see that $\rho(X_1),\ldots,\rho(X_n)$ are commuting matrices whose joint spectrum encodes Young tableaux, which in turn encodes the multiplicities of the irreducible representations.  In particular, Young tableaux appear in the spectrum as {\it content vectors}:


\Definition {Cont} A {\it content vector}  $\vc \alpha = \vect[a_1,\ldots,a_n] \in \Cont(n) \subset \bbZ^n$ satisfies the following:

\beginorderedlist
\item $a_1= 0$.
\item $\set{a_k - 1, a_k+1} \cap \set{a_1,\ldots,a_{k-1}} \ne \emptyset$ for $k = 2,\ldots,n$.
\item If $a_k = a_j = a$ for some $ k < j$ then
$$
\set{a-1,a+1} \subset \set{a_{k+1},\ldots,a_{j-1}}.
$$
\endorderedlist

There is a bijection between content vectors and Young Tableaux, \cite[Proposition 5.3]{Okounkov},  which we denote $\cont : \Y_n \rightarrow \Cont(n)$. This map can be constructed as follows: if the box in the $i$-th row and $k$-th column in a Young tableau has the number $k$ in it then the $k$-th entry of the corresponding content vector is $j-i$.  For example, for the Young diagram with partition $(3,2,1,1)$ the contents are
$$
\begin{ytableau}
0 & 1 & 2 \\
-1 & 0 \\
-2  \\
-3
\end{ytableau}
$$
and so the Young tableau in  \eqref{eq:ytableauexample}  corresponds to the content vector $\vect[0,1,-1,-2,0,2,-3]$. The map from content vectors to Young tableaux consists of filling in the boxes in the order the diagonals appear. In particular, if in the $k$th entry of the content vector we have an integer $p$ which has appeared $\mu$ times then, if $p$ is non-negative, the $(\mu+1,\mu+p+1)$ box is equal to $k$, otherwise, if $p$ is negative, the $(\mu-p+1,\mu+1)$ box is equal to $k$.

We now come to the key result: $\rho_\lambda(X_j)$ are diagonal for all irreducible representations and their entries encode all content vectors associated with the corresponding partition. The following is equivalent to \cite[Theorem 5.8]{Okounkov}:

\Lemma{diag}  $\rho_\lambda(X_j)$ is diagonal and {\addition each position on the main diagonal gives a content vector}: for each $k = 1,\ldots,d_\lambda$
$$
\Vectt[{\bfe_k^\top \rho_\lambda(X_1) \bfe_k}, \dots, {\bfe_k^\top \rho_\lambda(X_n) \bfe_k}] = \cont(y_k)
$$
where $\set{y_1, \ldots, y_{d_\lambda}} = \Y_n$ are an enumeration of  Young tableaux. 


%
%
%
%
%
%

\Section{rep2la} From representation theory to linear algebra.

%
%
%
%

A classic result in representation theory is that all representations of finite groups are decomposable:

\Proposition{directsum}\cite[Proposition 1.8]{FultonRepTheory}
For any representation $\rho : G \rightarrow \GL(V)$ of a finite group $G$ there is a decomposition
$$
	V \cong V_1^{\oplus a_1} \oplus \cdots \oplus V_r^{\oplus a_r}
$$
where  $\rho_k : G \rightarrow \GL(V_k)$ are distinct irreducible representations. The decompositions of $V$ into a direct sum of the $r$ summands is unique, as are the $V_k$ that occur (up to isomorphism) and their multiplicities $a_k$.


In terms of concrete linear algebra and specialising to the case of $G = S_n$ where we are given an orthogonal matrix representation $\rho : S_n \rightarrow \Orth(d)$, we denote the $a_k$ subspaces $V_k^{(1)}, \ldots, V_k^{(a_k)} \subset \bbR^d$ whose intersection is $\{0\}$ that are all isomorphic to  $V_k = \bbR^{d_k}$ where $d_k := d_{\lambda^k}$ for $\lambda^k \vdash n$.

We shall use Schur's Lemma to show that $\{B_k^{(p)}\}$  have certain orthogonality properties:

\Lemma{Schur}
\cite[Schur's Lemma 1.7]{FultonRepTheory} If $\rho_V : G \rightarrow \GL(V)$ and $\rho_W  : G \rightarrow \GL(W)$ are irreducible representations of $G$ and $\varphi : V \rightarrow W$ is a $G$-module homomorphism, then
\beginorderedlist
\item Either $\varphi$ is an isomorphism, or $\varphi = 0$.
\item If $\rho_V = \rho_W$, then $\varphi = cI$ for some $c \in  \bbC$.
\endorderedlist

This lemma encodes the fact that bases corresponding to differing irreducible representations are automatically orthogonal to each other whereas bases corresponding to the same irreducible representation have inner products that are a scaled identity.   {\addition We can use this to show that we can choose a basis for each $V_k^{(p)}$ so that the resulting representation is precisely of the form of \defref{irreps}:}

\Proposition{orthoBk} {\addition There exists $B_k^{(p)} \in  \bbR^{d \times d_k}$ with orthogonal columns, i.e.,  $(B_k^{(p)})^\top B_k^{(p)} = I_{d_k}$, such that $V_k^{(p)} = \colspan B_k^{(p)}$ and for $\rho_k := \rho_{\lambda^k}$ we have
$$
\rho(g) B_k^{(p)} = B_k^{(p)} \rho_k(g).
$$
}

{\addition

\Proof

Consider an orthogonal basis which spans $V_k^{(p)}$: if $V_k$ has dimension $d_k$ then we can construct a matrix $\tilde B_k^{(p)} \in  \bbR^{d \times d_k}$  with orthogonal columns that span $V_k^{(p)}$ and note that invariance guarantees that there exists a representation $\tilde \rho_k : S_n \rightarrow  \bbR^{d_k \times d_k}$ such that
$$
\rho(g) \tilde  B_k^{(p)} = \tilde  B_k^{(p)} \tilde  \rho_k(g).
$$
Because $V_k^{(p)}$ are isomorphic to $V_k$, there exists a matrix $U_k \in \bbR^{d_k \times d}$ that is invertible on $V_k^{(p)}$ so that for all $\bfc \in V_k^{(p)}$ we have $U_k  \rho(g) \bfc =  \rho_k(g) U_k \bfc$ therefore consider
$$
C :=  \tilde B_k^{(p)} (U_k  \tilde B_k^{(p)})^{-1}.
$$
This satisfies
$$
C \rho_k(g) = \tilde B_k^{(p)} (\rho_k(g^{-1}) U_k  \tilde B_k^{(p)})^{-1} = 
 \tilde B_k^{(p)} ( U_k \rho(g^{-1})  \tilde B_k^{(p)})^{-1}  = 
 \tilde B_k^{(p)} ( U_k   \tilde B_k^{(p)} \tilde \rho_k(g^{-1}))^{-1} = \rho(g)  C.
$$
Orthogonality follows from Schur's lemma. In particular, consider the map $\varphi : \bbR^{d_k} \rightarrow \bbR^{d_k}$ defined by $\varphi(v) :=  C^\top  C v$. This is a $G$-module homomorphism:
$$
\rho_k(g) \varphi(v) = (C \rho_k(g^{-1}))^\top C v = (\rho(g^{-1}) C)^\top C v = C^\top C \rho_k(g) v = \varphi(\rho_k(g) v).
$$
Thus $C^\top C = c I$ for some real constant $c$. Thus define
$$
B_k^{(p)} := C/\sqrt{c}.
$$
\mqed

The orthogonality property carries over to bases corresponding to different irreducible representations:

}
%

\Corollary{orthobasis} 
$$
(B_k^{(p)})^\top B_j^{(q)} = \begin{cases}
			 0_{d_k \times d_j} & k \ne j \\
			 m_k^{(pq)} I_{d_k} & k = j
			 \end{cases}
$$
for some constants $m_k^{(pq)} \in \bbR$.

\Proof
Consider the map $\varphi : \bbR^{d_j} \rightarrow \bbR^{d_k}$ defined by $\varphi(v) :=  (B_k^{(p)})^\top  B_j^{(q)} v$. {\addition Similar to above this is also} a $G$-module homomorphism:
$$
\rho_k(g) \varphi(v) = (B_k^{(p)} \rho_k(g^{-1}))^\top B_j^{(q)} v = (\rho(g^{-1}) B_k^{(p)})^\top B_j^{(q)} v = (B_k^{(p)})^\top B_j^{(q)} \rho_j(g) v = \varphi(\rho_j(g) v).
$$
The corollary therefore follows from Schur's lemma, where the fact that $m_k^{(pq)}$ is real follows since $B_k^{(p)}$ and $B_j^{(q)}$ have real entries.
\mqed

The previous proposition guarantees that if all irreducible representations have trivial multiplicity ($a_1 = \cdots = a_r = 1$) then $U = \bvect[B_1^{(1)}, B_2^{(1)}, \ldots, B_r^{(1)}]$ is an orthogonal matrix. This is not necessarily the case when we have non-trivial multiplicities. However, we can guarantee the existence of an orthogonal matrix that fully block-diagonalises a representation by taking an appropriate linear combination of $B_k^{(1)},\ldots,B_k^{(a_k)}$, thus guaranteeing that \probref{blockdiag} has a solution:

\Theorem{fullblockdiag} There exists an orthogonal matrix $Q$ such that
$$
Q^\top \rho(g) Q = \sopmatrix{ \rho_1(g)^{\oplus a_1} \\ & \ddots \\ && \rho_r(g)^{\oplus a_r}}
$$
where $\rho_1,\ldots,\rho_r$ are all irreducible representations. 

\Proof

 Recall $(B_k^{(p)})^\top B_k^{(q)} = m_k^{(pq)} I_{d_k}$ and consider the Cholesky decomposition
$$
\underbrace{\begin{pmatrix}
m_k^{(11)} & \cdots & m_k^{(1a_k)} \\
\vdots & \ddots & \vdots \\
m_k^{(a_k1)} & \cdots & m_k^{(a_ka_k)} 
\end{pmatrix}}_{M_k} = \tilde R_k^\top  \tilde R_k
$$
{\addition where $\tilde R_k$ is upper-triangular}. Note that $\tilde R_k$ is invertible as $M_k$ is a Gram matrix (associated with the first columns of $B_k^{(1)},\ldots,B_k^{(a_k)}$). 
Define $R_k = \tilde R_k \otimes I_{d_k}$ and  $Q_k := \bvect[B_k^{(1)},\ldots,B_k^{(a_k)}] R_k^{-1}$, which has orthogonal columns as $Q_k^\top Q_k = R_k^{-\top} (M_k \otimes I_{d_k}) R_k^{-1}  = I_{a_k d_k}$.
We then have
\meeq{
\rho(g) Q_k = \bvect[B_k^{(1)},\ldots,B_k^{(a_k)}]  (I_{a_k} \otimes \rho_k(g)) (\tilde R_k^{-1} \otimes I_{d_k}) \ccr
= \bvect[B_k^{(1)},\ldots,B_k^{(a_k)}]   (\tilde R_k^{-1} \otimes I_{d_k})  (I_{a_k} \otimes \rho_k(g)) 
= Q_k   \rho_k(g)^{\oplus a_k}.
}
Thus $Q := \bvect[Q_1, \ldots, Q_r]$ satisfies the necessary properties. 
\mqed

\Section{mult} Counting multiplicities via the Gelfand--Tsetlin algebra and basis.

We now consider the solution of \probref{multiplicities}: how do we use the existence of an orthogonal basis to determine how many copies of each irreducible representation are present in a given representation? The key will be the simultaneous diagonalisation of related commuting operators: the representations of the YJM-generators $\rho(X_1),\ldots,\rho(X_n)$.

Let us return to the {\addition example representation $\rhoe : S_4 \rightarrow O(6)$}, where we have: {\addition
\meeq{
\rhoe(X_1) = 0,  
&\rhoe(X_2) & =   \sopmatrix{ 0 & 1 \cr 1 & 0 \cr && 1 \cr &&& 1 \cr &&&& 0 & 1 \\ &&&& 1 & 0 },  \ccr
\rhoe(X_3)  = \sopmatrix{ 1 & 0 & 1 \cr 0 & 1 & 1 \cr 1 & 1 & 0 \cr &&& 2 \\ &&&& 0 & 2 \\ &&&& 2 & 0  }, 
&\rhoe(X_4)  &= \sopmatrix{ 2 & 0 & 0 & 1 \cr 0 & 2 & 0 & 1 \cr 0 & 0 &  2 & 1  \cr 1 & 1&1 & 0 \\
&&&& 0 & 3\\
&&&& 3 & 0
}.
}
}
Note the above are all symmetric matrices. In fact, they also commute. If we conjugate with the orthogonal matrix $Q$ from \eqref{eq:permutationQ} that block-diagonalised the representation we have that the GZ generators are simultaneously diagonalised: {\addition
\meeq{
Q^\top \rhoe(X_1) Q = \mdiag(0,0,0,0,0,0), \ccr
Q^\top \rhoe(X_2) Q = \mdiag(-1,-1,1,1,1,1), \ccr
Q^\top \rhoe(X_3) Q = \mdiag(-2,1,-1,2,2,2),\ccr
 Q^\top \rhoe(X_4) Q =\mdiag(-3,2,2,-1,3,3). 
}
}
This property is general:

\Proposition{comm} If $\rho$ is an orthogonal matrix representation then $\rho(X_1),\ldots,\rho(X_n)$ are symmetric and commute. In other words, there exists an orthogonal matrix  that simultaneously diagonalises $\rho(X_1), \ldots, \rho(X_n)$. 

\Proof

One approach to showing this is to note that $(i\ j) = \tau_i \tau_{i+1} \cdots \tau_{j-1} \tau_j \tau_{j-1} \cdots \tau_i$ hence $\rho((i\ j))$ is a symmetric matrix, and therefore so is $\rho(X_j)$. 
The fact that $\rho(X_j)$ commute follows since  $X_j$ commute as outlined in \cite{Okounkov}. 

Another approach that is more illuminating to our problem is to appeal to \thref{fullblockdiag}: using the $Q$ that fully block-diagonalises $\rho(g)$ we have
\meeq{
Q^\top \rho(X_j) Q =  Q^\top \rho((1\ j+1)) Q  + \cdots + Q^\top \rho((j\ j+1)) Q \ccr
=  \sum_{\ell=1}^j \sopmatrix{ \rho_1((\ell\ j+1))^{\oplus a_1} \\ & \ddots \\ && \rho_r((\ell\ j+1))^{\oplus a_r}} \ccr
= \sopmatrix{ \rho_1(X_j)^{\oplus a_1} \\ & \ddots \\ && \rho_r(X_j)^{\oplus a_r}} 
}
and \lmref{diag} guarantees that $\rho_k(X_j)$ are diagonal.
\mqed

We can use the fact that the eigenvalues of $\rho(X_j)$ contain $a_k$ copies of the diagonal entries of $\rho_k(X_j)$ to deduce the multiplicities $a_k$ from the joint spectrum.

\Definition{Spec}  Given an orthogonal matrix representation $\rho$,
 $\Spec_\rho \in \bbR^{d \times n}$ is the matrix containing the joint spectrum of $\rho(X_1), \ldots, \rho(X_n)$. That is, if $Q = \bvect[\bfq_1, \ldots, \bfq_d] \in \Orth(d)$ simultaneously diagonalises $\rho(X_1),\ldots,\rho(X_n)$ we have
 $$
\vc e_k^\top \Spec_\rho \vc e_j = \bfq_k^\top \rho(X_j) \bfq_k.
 $$
Note this is only uniquely defined up to permutation of rows: we choose the ordering so that the first column is non-decreasing, rows with the same entry in the first column are non-decreasing in the second column, and so-forth.

{\addition

We can deduce which irreducible representations are present from $\Spec_\rho$:

\Corollary{rowsarecont} The rows of $\Spec_\rho$ are content vectors corresponding to each Young tableau associated with the irreducible representation $\rho_k$, repeated $a_k$ times.

\Proof

This follows from the (second) proof of \propref{comm}.

\mqed
}

In the case of the {\addition example} representation we rearrange the eigenvalues from above into a matrix:{\addition
$$
\Spec_{\rhoe} = \sopmatrix{
0 & -1 & -2 & -3 \\
0 & -1 & 1 & 2 \\
0 & 1 & -1 & 2 \\
0 & 1 & 2 & -1 \\
0 & 1 & 2 & 3 \\
0 & 1 & 2 & 3
}
$$
}

The rows of $\Spec_\rho$ are content vectors, whose corresponding Young tableaux according to the map $\cont^{-1}$ are:
$$
{\addition \begin{ytableau}
    1 \\
    2 \\ 
    3 \\ 
    4
    \end{ytableau} },
\begin{ytableau}
    1 & 3 & 4 \\
    2
    \end{ytableau}, \begin{ytableau}
    1 & 2 & 4 \\
    3
    \end{ytableau},\begin{ytableau}
    1 & 2 & 3 \\
    4
    \end{ytableau}, \begin{ytableau}
    1 & 2 & 3 & 4
    \end{ytableau}, {\addition \begin{ytableau}
    1 & 2 & 3 & 4
    \end{ytableau}}.
$$
That is, we have every possible Young tableau corresponding to the partitions {\addition $\lambda^1 = 1 + 1 + 1 + 1$, $\lambda^2 = 3 + 1$ and $\lambda^3 = 4$ (which is repeated twice)}. We can deduce from this that each of the corresponding multiplicities {\addition are $a_1 = a_2 = 1$ and $a_3 = 2$}.

We are yet to discuss how to simultaneously diagonalise $\rho(X_1),\ldots,\rho(X_n)$. Practically speaking, this is a well-studied problem with effective
iterative methods \cite{bunse1993numerical}. However, we can use the fact that the eigenvalues are integers to ensure that the problem is solvable with a finite number of operations via a more traditional {\it Diagonalize-One-then-Diagonalize-the-Other (DODO)} approach. 

\Lemma{cubicjointspec} 
$\rho(X_1),\ldots,\rho(X_n)$ can be simultaneously diagonalised in $\Ord(n^2 d^3)$ operations.

\Proof

We first note that the eigenvalues of $\rho(X_j)$ all lie between $1-j$ and $j-1$ {\addition since they correspond to entries of content vectors}. A symmetric matrix $A$ with integer eigenvalues satisfying $1-j \leq \lambda \leq 1-j$ can be diagonalised in $\Ord(d^3 + jd)$ operations: first one can tridiagonalize
$$
Q_H^\top A Q_H = T
$$ 
where $Q_H$ is a product of Householder reflections and $T$ is a symmetric tridiagonal matrix with the same eigenvalues as $A$
in $\Ord(d^3)$ operations. Determination of an orthogonal basis for the nullspace of a tridiagonal matrix requires $\Ord(d)$ operations using e.g. Gram--Schmidt, thus $T$ can be diagonalised in $\Ord(j d)$ operations by determination of the nullspaces of $T - m I$ for each $m = 1-j, \ldots, j-1$.
Thus we can diagonalise $\rho(X_2) = Q_2 \Lambda_2 Q_2^\top$ in $\Ord(d^3)$ operations (we know it has eigenvalues $\pm1$) and conjugate $Q_2^\top \rho(X_j)Q_2$ for $j  = 3,\ldots,n$ in a total of $\Ord(n d^3)$ operations. These will be block-diagonalised according to the eigenvalues of $\rho(X_2)$ hence we can deflate these matrices and repeat the process $n$ times on the sub-matrices.
\mqed

By converting the content vectors associated with the joint spectrum to partitions (see \algref{Cont2Part}) we can deduce the multiplicities of the irreducible representations:

\Corollary{prob1} \probref{multiplicities} can be solved in $\Ord(n^2 d^3)$ operations.

{\addition \Remark In practice, this complexity can be improved using randomised linear algebra, eg., computing a single eigenvalue decomposition of a randomised linear combination of $X_j$ \cite{he2024randomized} would reduce the complexity to $\Ord(n d^3)$.

}

\Section{blockdiag} Decomposing multiple copies of the same representation.

The above results guarantee that one can compute a $\tilde Q$ that simultaneously diagonalises the symmetric matrices $\rho(X_1),\ldots,\rho(X_n)$. Moreover, we have block-diagonalised $\rho$, but unfortunately multiple copies of the same irreducible representation are not necessarily decoupled. That is, $\tilde Q$ is only guaranteed to reduce a representation to a {\it canonical representation}:

\Lemma{cano} Suppose $\tilde Q$  simultaneously diagonalises the YJM generators $\tilde Q^\top \rho(X_j) \tilde Q = \Lambda_j$, where we assume the $\Lambda_j$ are sorted by the corresponding partitions.  Then it reduces $\rho(g)$ to a canonical representation:
$$
 \rho(g) \tilde Q =  \tilde Q \sopmatrix{ \tilde \rho_1(g) \\ & \ddots \\ && \tilde \rho_r(g)}
$$
where there exists  matrices $\bar Q_k \in \Orth(a_k d_k)$ for irreducible representation $\rho_k$ of dimension $d_k$ such that
$$
\tilde \rho_k(g) \bar Q_k  = \bar Q_k \rho_k(g)^{\oplus a_k}.
$$
Moreover, if $\bar Q_k$ is divided into $a_k$ blocks of size $d_k \times d_k$ then each block is diagonal. 

\Proof

We know that $Q = \bvect[Q_1,\dots,Q_r]$ from \thref{fullblockdiag} and $\tilde Q = \bvect[\tilde Q_1,\dots,\tilde Q_r]$ are both orthogonal matrices that simultaneously diagonalise $\rho(X_j)$, so the only ambiguity arrives due to each content vector being repeated according to the multiplicity of the irreducible representation. Thus the columns of $\tilde Q_k$ associated with a given content vector must be linear combinations of the columns of $Q_k$ associated with the same content vector, which ensures that $\bar Q_k = \tilde  Q_k^\top Q_k$ has diagonal blocks.
\mqed

Thus we further need a method to decompose a representation containing multiple copies of the same irreducible representation. As we have reduced the representation to a canonical representation we can consider each $\tilde \rho_k$ separately, which for simplicity we denote $\tilde \rho$ with corresponding partition $\lambda$.   We mimic an approach for building an orthogonal basis for the eigenspace of a matrix whose corresponding eigenvalue is non-trivial: that is, the problem of finding  orthogonal vectors ${\bf q}_k$ such that $(A - \lambda I){\bf q}_k = 0$. We instead find orthogonal matrices $Q^{(k)}$ such that $\tilde \rho(g) Q^{(k)} - Q^{(k)} \rho_\lambda(g) = 0$. Both problems  consist of finding the nullspace of a matrix. A remarkable fact is that solving this nullspace problem in the na\"ive way guarantees orthogonality.

\Theorem{multnullspace} 
Given an orthogonal matrix representation $\tilde \rho : S_n \rightarrow \Orth(a d_\lambda)$, we have
\begin{equation}\label{eq:multrho}
\tilde \rho(g) \bar Q = \bar Q  \rho_\lambda(g)^{\otimes a}
\end{equation}
for an orthogonal matrix $\bar Q = \bvect[\bar Q^{(1)},\ldots,\bar Q^{(a)}]$ with $\bar Q^{(\ell)} \in \bbR^{a d_\lambda \times d_\lambda}$ if and only if 
\begin{equation}\label{eq:kronQ}
\underbrace{\begin{pmatrix}
\tilde \rho(\tau_1) \otimes I_{d_\lambda}- I_d \otimes \rho_\lambda(\tau_1) \\
\vdots \\
\tilde \rho(\tau_{n-1}) \otimes I_{d_\lambda}- I_d \otimes \rho_\lambda(\tau_{n-1})
\end{pmatrix}}_{\addition K_{\tilde \rho, \lambda}} U = 0
\end{equation}
for a matrix $U = \bvect[\bfu^{(1)}, \ldots, \bfu^{(a)}] \in \bbR^{(ad_\lambda^2) \times a}$ with orthogonal columns satisfying
$$
 \vec \bar Q^{(\ell)} = \sqrt{d_\lambda}\bfu^{(\ell)},
$$
where $\vec : \bbR^{m \times n} \rightarrow \bbR^{mn}$ concatenates the columns of a matrix.

\Proof

Given $\bar Q$, define  $U = \bvect[\bfu^{(1)}, \ldots, \bfu^{(a)}]$ with $\bfu^{(\ell)} = {\vec \bar Q^{(\ell)} / \sqrt{d_\lambda} }$. The fact that it satisfies \eqref{eq:kronQ} follows immediately from the definition of the Kronecker product so we only need to show orthogonality. Writing $\bar Q^{(\ell)} = \bvect[\bfq_1^{(\ell)},\ldots,\bfq_{d_\lambda}^{(\ell)}]$ we have for $k \ne \ell$
$$
(\bfu^{(k)})^\top \bfu^{(\ell)} = {(\bfq_1^{(k)})^\top \bfq_1^{(\ell)} + \cdots + (\bfq_{d_\lambda}^{(k)})^\top \bfq_{d_\lambda}^{(\ell)} \over d_\lambda} = 0
$$
and (because the Frobenius norm squared of a matrix with orthogonal columns is equal to its rank)
$$
\| \bfu^{(k)} \| = {1 \over \sqrt d_\lambda} \| \bar Q^{(k)} \|_F = 1.
$$

For the other direction, define $\bar Q^{(\ell)} \in \bbR^{ad_\lambda \times d_\lambda}$ by $\vec \bar Q^{(\ell)} = \sqrt{d_\lambda} \bfu^{(\ell)}$ (that is, we reshape the vector $\bfu^{(\ell)}$ to be a matrix) and define $\bar Q = \bvect[\bar Q^{(1)},\ldots,\bar Q^{(a)}]$. Note that \eqref{eq:multrho} is automatically satisfied by the definition of the Kronecker product so we need only show orthogonality. 
From Schur's lemma (\lmref{Schur}) we know that $(\bar Q^{(k)})^\top \bar Q^{(\ell)} = c_{k\ell} I $ for some constants $c_{k\ell}$. Writing $\bar Q^{(\ell)} =  \bvect[\bfq_1^{(\ell)},\ldots,\bfq_{d_\lambda}^{(\ell)}]$ this means that $(\bfq_s^{(k)})^\top \bfq_t^{(\ell)} =   c_{k\ell}  \delta_{st}$. Thus we have
$$
\delta_{k\ell} = (\bfu^{(k)})^\top \bfu^{(\ell)} ={\addition {\sum_{j = 1}^{d_\lambda} (\bfq_j^{(k)})^\top \bfq_j^{(\ell)}  \over  d_\lambda} =  {\sum_{j = 1}^{d_\lambda} c_{k\ell} \over  d_\lambda} } =  c_{k\ell}
$$
which  ensures that $\bar Q$ is in fact orthogonal.
\mqed

%

{\addition 
\Corollary{subcols} 
In the case where $\tilde \rho_k$ is defined as in \lmref{cano} the nullspace of $K_{\tilde \rho_k, \lambda_k}$ as defined in \eqref{eq:kronQ} can be recovered from the nullspace of only $a_k d_k$ of its columns. 

\Proof

From \lmref{cano} we know that $\bar Q_k = \bvect[\bar Q_k^{(1)}, \dots, \bar Q_k^{(a_k)}]$ has diagonal blocks and hence has only $a_k^2 d_k$ non-zero entries, or in particular each $\bar Q_k^{(\ell)}$ has $a_k d_k$ non-zero entries, in the same locations. Dropping the known zero rows of $\vec \bar Q_k^{(\ell)}$ will therefore be in the nullspace of $K_{\tilde \rho_k, \rho_k}$ with the corresponding columns dropped.

\mqed

}

\Corollary{fullyreduce} 
	\probref{blockdiag} can be solved in $\Ord(n d^4)$ operations.

\Proof
Note that	$\bar Q_k^{(1)}, \ldots, \bar Q_k^{(a_k)}$ can na\"\i vely be computed in $\Ord(n a_k^3 d_k^5)$ operations via Gram--Schmidt applied to the $\bbR^{n a_k d_k^2 \times a_k d_k^2}$ matrix in \eqref{eq:kronQ}, for a total of $\Ord(n d^5)$ operations (using $\sum_{k=1}^r a_k d_k = d$).  {\addition Using the previous corollary allows us to reduce the complexity by applying Gram--Schmidt to an $\bbR^{n a_k d_k^2\times a_k d_k}$ matrix, taking  $\Ord(n  a_k^3 d_k^4)$ operations.} Summing over these for $k=1,\ldots,r$ completes the proof.
\mqed

{\addition 
\Remark In practice many of the rows of the matrix $K_{\tilde \rho_k, \lambda_k}$ are identically zero, and can be dropped without altering its nullspace. Such rows can be found in $\Ord(n a_k^2 d_k^3)$ operations. In experiments this appears to reduce the overall complexity to $\Ord(n d^3)$, see \secref{examples}, however, if this is true in general remains open and would require understanding the sparsity present in the generators of the irreducible representations $\rho_{\lambda_k}$. 

}

\Section{algorithm} Algorithms. 

Encoded in the above results  are algorithms for solving Problems 1--3.  Here we outline explicitly the stages of the algorithms and discuss briefly the practical implementation in floating point arithmetic. \algref{multiplicities} computes irreducible representation multiplicities, that is, it solves \probsref[present,multiplicities], as well as computing an orthogonal matrix $\tilde Q$ that reduces a representation to canonical representations. \algref{blockdiag} builds on \algref{multiplicities} to fully decompose a representation into irreducible representations, solving \probref{blockdiag}, including the case where there are non-trivial multiplicities. 
{\addition
For completeness we also include \algref{YJM}, which constructs the YJM-generators, and  \algref{Cont2Part}, which discusses the translation of a content vector to a partition. 
}

\begin{algorithm}[tb]
\caption{Compute multiplicities of irreducible representations (\probref{multiplicities}) \label{Algorithm:multiplicities}}
{\noindent\bf Input:} {Generators $\rho(\tau_1), \ldots, \rho(\tau_{n-1}) \in \Orth(d)$.}

{\noindent\bf Output:} {Partitions $\lambda^1,\ldots,\lambda^r$ where $\lambda^k \vdash n$, corresponding multiplicities $a_1,\ldots,a_r$, and an orthogonal matrix $\tilde Q \in \Orth(d)$ that reduces a representation to a canonical representation.}

\begin{algorithmic}[1]
\State Construct $\rho(X_1)=0, \rho(X_2),\ldots,\rho(X_n) \in \bbR^{d \times d}$ {\addition using \algref{YJM}}.
\State Find $Q$ that simultaneously diagonalises $\rho(X_2),\ldots,\rho(X_n)$, that is $\tilde Q^\top \rho(X_k) \tilde Q = \Lambda_k$ where $\Lambda_k$ is diagonal. This can be accomplished via procedure outlined in the proof of \lmref{cubicjointspec}.
\State The diagonals of $\Lambda_1 = 0, \Lambda_2, \ldots, \Lambda_n$ {\addition give $\Spec_\rho$ whose rows are} are content vectors. Convert the content vectors to partitions {\addition using \algref{Cont2Part}} to determine $\lambda^1, \ldots, \lambda^r$ which are repeated $\tilde a_1, \ldots, \tilde a_r$ times. Each partition will be repeated once for each corresponding Young tableau so the  multiplicities are $a_k = \tilde a_k/d_{\lambda^k}$. 
\end{algorithmic}

\end{algorithm}

\begin{algorithm}[tb]
\caption{Decompose representations into irreducible representations (\probref{blockdiag})}\label{Algorithm:blockdiag}

{\noindent\bf Input:} {Generators $\rho(\tau_1), \ldots, \rho(\tau_{n-1}) \in \bbR^{d \times d}$.}

{\noindent\bf Output:} {An orthogonal matrix $Q \in \Orth(d)$ that fully decomposes a representation.}

\begin{algorithmic}[1]
\State Find $\lambda^1,\ldots,\lambda^r$,  $a_1,\ldots,a_r$, and $\tilde Q = \bvect[\tilde Q_1,\dots, \tilde Q_r]$ using \algref{multiplicities}.
\For{$k = 1, \ldots, r$}
\State For the representation $\tilde \rho_k(g) = \tilde Q_k^\top \rho(g) \tilde Q_k$, solve  \eqref{eq:kronQ} using Gram–Schmidt to determine the matrix $U_k = \bvect[\bfu_k^{(1)}, \ldots, \bfu_k^{(a_k)}]$. {\addition In order to reduce the complexity of  Gram--Schmidt  we drop the columns of $K_{\tilde \rho_k,\lambda_k}$ corresponding to guaranteed zero entries according to \corref{subcols} and drop any rows detected to be zero.}
\State Construct $\bar Q_k^{(k)} \in \bbR^{d \times d_\lambda}$ such that  $\vec \bar Q_k^{(k)} = \sqrt{d_\lambda}\bfu_k^{(k)}$ (that is, scale and reshape the vector $\bfu_k^{(k)}$). Then $Q_k = \tilde Q_k \bvect[\bar Q_k^{(1)}, \dots, \bar Q_k^{(a_k)}]$.
\EndFor
\State  $Q = \bvect[Q_1, \dots, Q_r]$.
\end{algorithmic}
\end{algorithm}

\begin{algorithm}[tb]

\caption{\addition Compute representations of the Young--Jucys--Murphy (YJM) generators \label{Algorithm:YJM}}
\addition
{\noindent\bf Input:} {Generators $\rho(\tau_1), \ldots, \rho(\tau_{n-1}) \in \Orth(d)$.}

{\noindent\bf Output:} {YJM generators  $\rho(X_2), \ldots, \rho(X_n)$.}

\begin{algorithmic}[1]
\For{$j = 2, \ldots, n$} 
\State Let $b_{j,j-1} = a_{j,j-1} = \rho(\tau_{j-1})$.
	\For{$k = j-2, \ldots, 1$}
	\State Let $a_{k,j} = \rho(\tau_k) a_{k+1,j} \rho (\tau_k)$ and $b_{k,j} = b_{k+1,j} + a_{k,j}$.
	\EndFor
	\State  $\rho(X_j) = b_{1,j}$	.
\EndFor
\end{algorithmic}
\end{algorithm}

\begin{algorithm}[tb]
\caption{\addition Convert a content vector to a partition \label{Algorithm:Cont2Part}}
\addition
{\noindent\bf Input:} {Content vector $c \in \Cont(n)$}

{\noindent\bf Output:} {Partition $\lambda \vdash n$}

\begin{algorithmic}[1]
\State  The partition has $m = 1-\min c$ components. Initialise $\lambda_1 = \ldots = \lambda_m = 0$. 
\For{$k = 0,-1,-2,\ldots,\min c$}
	\State If $k$ occurs $\mu$ times in $c$, increment $\lambda_{1-k}, \ldots, \lambda_{\mu-k}$ by one.
\EndFor	
\For{$k = 1,2,\ldots,\max c$}
	\State  If $k$ occurs $\mu$ times in $c$, increment $\lambda_1, \ldots, \lambda_\mu$ by one.	
\EndFor

\end{algorithmic}
\end{algorithm}

In the practical implementation \cite{NumericalRepresentationTheory} we use floating point arithmetic, and in particular we simultaneously diagonalise matrices using standard numerical  methods {\addition for diagonalising matrices}  (e.g. we could use \cite{bunse1993numerical} though in practice we use the simpler to implement {\it diagonalize-one-then-diagonalize-the-other (DODO)} approach), as opposed to the proposed exact method based on nullspace calculations. For the nullspace calculation in \algref{blockdiag} we use a standard method which is based on computing the Singular Value Decomposition (SVD) and taking the singular vectors associated with the {\addition smallest $a_k$ singular values, which will all be } approximately zero.   We also use sparse matrix data structures to further improve the computational cost. Using floating point arithmetic introduces round-off error, however, in practice the error is {\addition proportional to machine epsilon ($\approx  2.22 \times 10^{-16}$)} and the values of the approximate eigenvalues can be rounded to exact integers: that is, while there is no proof of correctness, when the algorithm succeeds it can be verified, in particular as we  have  computed  a $Q$ that approximately block-diagonalises the representation and we have precise bounds on the errors in floating point matrix multiplication \cite{higham2002accuracy}. Note however the dependence on black-box linear algebra software implemented in floating point means there is a small but non-zero chance of failure.

\Section{examples} Examples.

We now present some examples. In \figref{tensor_321_222} we apply  \algref{multiplicities} (using floating point arithmetic) to two representations resulting from tensor products: $\rho_{(3,2,1)} \otimes \rho_{(2,2,2)}$ (i.e. computing Kronecker coefficients)  and triple tensor product $\rho_{(3,2,1)} \otimes \rho_{(2,2,2)} \otimes \rho_{(3,3)}$. In \figref{kronpowtimes} we demonstrate the cubic complexity of applying \algref{multiplicities} as the dimension grows by considering increasing tensor powers
$$
\rho_\lambda^{\otimes k} = \underbrace{\rho_\lambda \otimes \cdots \otimes \rho_\lambda}_{\hbox{$k$ times}}.
$$
for fixed irreducible representation $\rho_\lambda$. 
Note for generic partitions $\lambda  \vdash n$ the dimension of the irreducible representation grows combinatorially fast as $n$ increases. However, for  {\it hook} irreducible representations,  which are those associated with partitions of the form $(m+1,1,1,\ldots,1)$, the growth is only quadratic. Thus in \figref{kronpowtimes} we focus on tensor powers of a hook and an almost-hook {\addition (i.e. one associated with the partition ${(2,2,1,1,\ldots,1,1)}$ where the number of ones are such that it is a partition of $n$)}  irreducible representations to scale to larger $n$. This example shows that the computational cost in practice primarily depends on the dimension of the representation, not $n$. 

\Figuretwoframed[tensor_321_222,tensor_321_222_33]
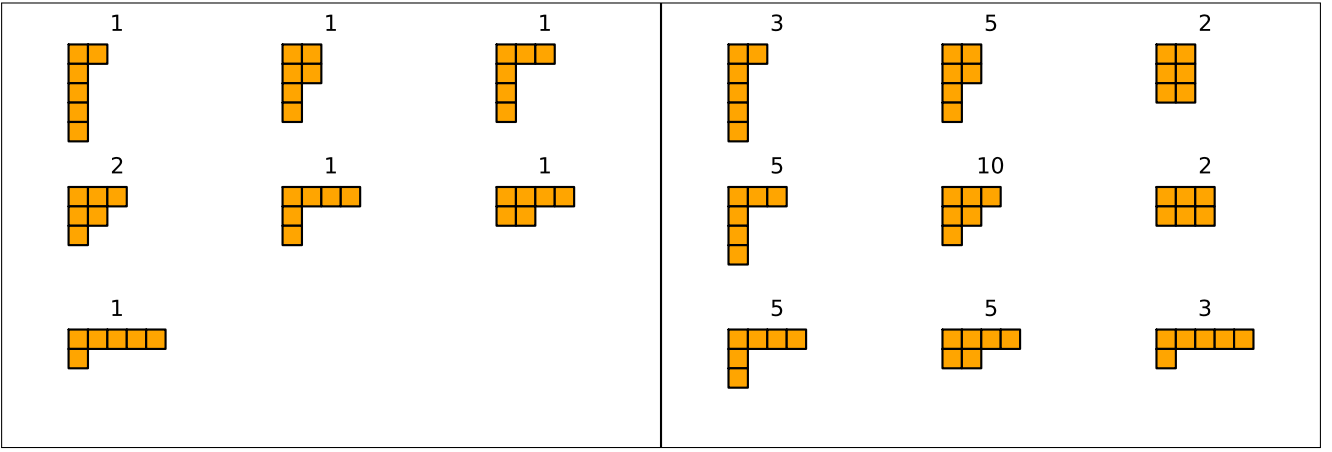
	Left: multiplicities of irreducible representations in the tensor product $\rho_{(3,2,1)} \otimes \rho_{(2,2,2)}$ (Kronecker coefficients). Right: multiplicities of irreducible representations in the triple tensor product $\rho_{(3,2,1)} \otimes \rho_{(2,2,2)} \otimes \rho_{(3,3)}$.
	

\Figurew{kronpowtimes}{0.45\hsize}
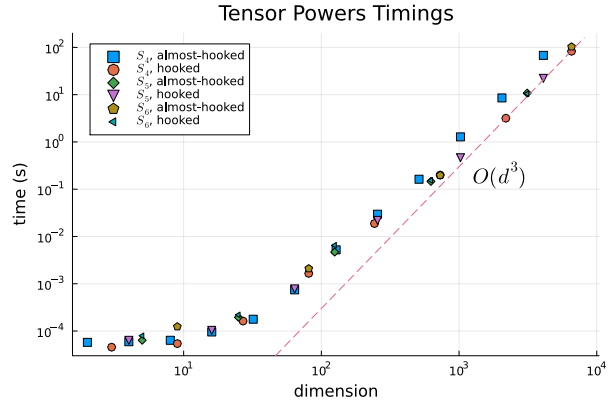
	Time taken to compute tensor powers ($\rho^{\otimes k}$) for  a hook ($\rho_{(2,1,1,\ldots,1,1)}$) and an almost-hook ($\rho_{(2,2,1,1,\ldots,1,1)}$) irreducible representation for varying $k$ and $n$. We plot the time taken compared to the dimension of the resulting representation to demonstrate the cubic growth in complexity.

For our next example, in \figref{cleschgordondimension} we consider the growth in computational cost  as $n$ increases for two examples: a hook tensor product $\rho_{(n-1,1)}  \otimes \rho_{(2,1,\ldots,1)}$ and an almost-hook tensor product $\rho_{(n-2,1,1)}  \otimes \rho_{(2,2,1,\ldots,1)}$. The dimensions of the irreducible representations $\rho_{(n-1,1)}, \rho_{(n-2,1,1)}$,  $\rho_{(2,1,\ldots,1)}$ {\addition and $\rho_{(2,2,1,\ldots,1)}$} are $n-1$, {\addition $(n-1)(n-2)/2$}, $n-1$ {\addition and $n(n-3)/2$}, respectively, so that the tensor products have dimensions  $(n-1)^2$ and {\addition $n(n-1)(n-2)(n-3)/4$}, hence we have demonstrated that the algorithm has  polynomial  complexity in $n$, in particular at worst $\Ord(n^7)$ {\addition and $\Ord(n^{13})$} operations.  In \figref{cleschgordondimension} we plot the time taken showing that the numerical implementation using floating point arithmetic {\addition demonstrating the predicted complexities}.

\Figuretwo[cleschgordondimension,cleschgordontiming]
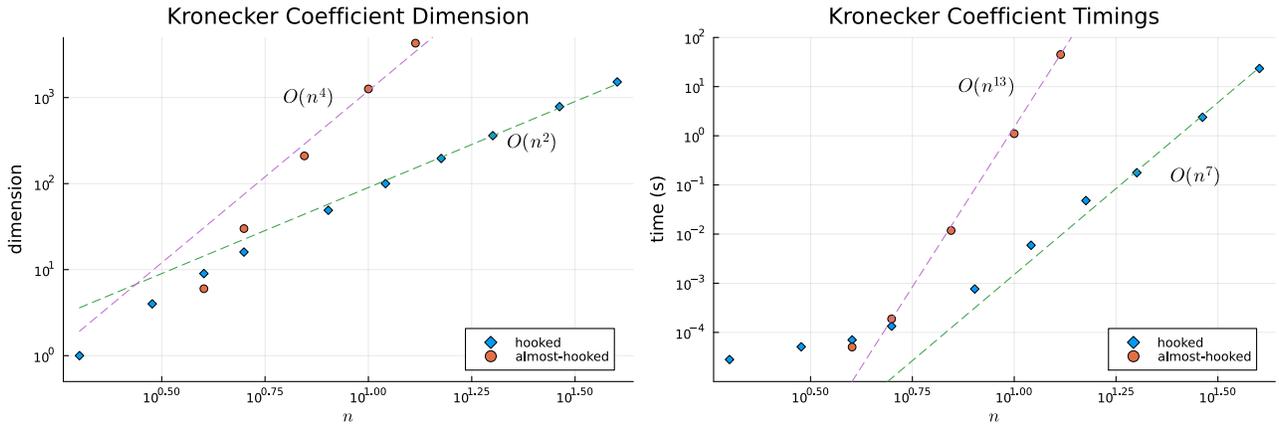
	Left:  growth in the dimension of  a hook  ($\rho_{(n-1,1)}  \otimes \rho_{(2,1,\ldots,1)}$) and an almost-hook ($\rho_{(n-2,1,1)}  \otimes \rho_{(2,2,1,\ldots,1)}$) tensor product. Right: time to compute corresponding Kronecker coefficients, {\addition showing the timings match the predicted rate of $\Ord(n^{7})$ and $\Ord(n^{13})$ operations}.

{\addition Finally, we consider an example that requires solving \probref{blockdiag} via \algref{blockdiag}. Let $p_k(x)$ denote a family of orthogonal polynomials\footnote{\addition What follows also applies to monomials, which can be viewed as orthogonal polynomials on the (complex) unit circle.} with respect to a weight $w(x)$ so that
$$
p_{{\bf k}}(\bfx) = p_{k_1}(x_1) \cdots p_{k_n}(x_n)
$$
are degree $k_1 + \cdots + k_n = p$ multivariate orthogonal polynomials with respect to a tensor product weight $w(x_1)\cdots w(x_n)$. Write all degree $p$ orthogonal polynomials in a vector ${\bfP}_p : \bbR^n \rightarrow \bbR^{d_{n,p}}$ where
$
d_{n,p} = \binom{p+n-1}{p} 
$
is the number of monomials of degree $p$ in $n$ variables. For any permutation $g$ we have that  ${\bfP}_p(\rhos_n(g) \bfx)$ is another basis of orthogonal polynomials with respect to the same weight.  As orthogonal polynomials are unique up to their span  \cite{dunkl2014orthogonal}, the basis must be an orthogonal linear combinations of the original basis and hence there exists matrices $\rho_{n,p}(g) \in O(d_{n,p})$ such that
$$
{\bfP}_p(\rhos_n(g) \bfx) = \rho_{n,p}(g)  {\bfP}_p(\bfx).
$$
In fact these matrices can be deduced directly as simple permutation matrices from the definition above.  

\Figuretwo[polyblockdiag_p,polyblockdiag_n]
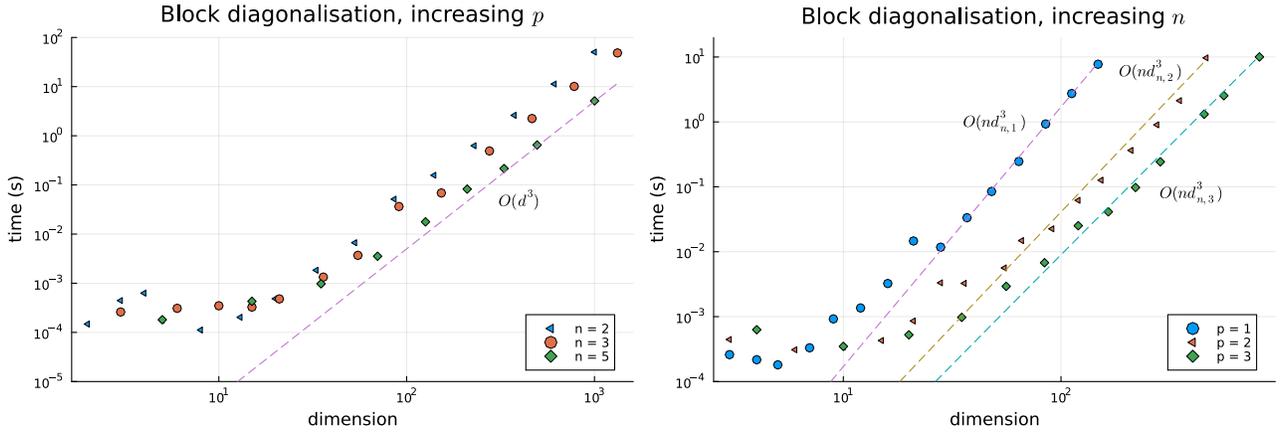
	{\addition Time taken to  block-diagonalise a representation generated from permuting variables of tensor product orthogonal polynomials in $n$ variables. The complexity of the implementation appears to be $\Ord(n d_{n,p}^3)$ operations, which is better than the proven complexity of $\Ord(n d_{n,p}^4)$, where $d_{n,p}$ is the dimension of the representation, in this case the number of monomials of degree $p$ in $n$ 
	variables.       Left:  we fix the number of variables $n$ and let the order of the polynomials increase. Right: we fix the order $p$ of the polynomials and increase the number of variables $n$}.

{\addition Note they have the requisite group structure:
\meeq{
\rho_{n,p}(g_2) \rho_{n,p}(g_1)  {\bfP}_p(\bfx) = \rho_{n,p}(g_2) {\bfP}_p(\rhos_n(g_1) \bfx) ={\bfP}_p( \rhos_n(g_2) \rhos_n(g_1) \bfx) = {\bfP}_p(\rhos_n(g_2 g_1) \bfx) \ccr 
=\rho_{n,p}(g_1 g_2) {\bfP}_p( \bfx)
}
thus $\rho_{n,p}$ defines a representation.  If we compute $Q \in O(d_{n,p})$ that block-diagonalise this representation using \algref{blockdiag} then ${\bfQ}_p(\bfx) :=  Q^\top {\bfP}_p(\bfx) $ is another basis of orthogonal polynomials but for one which permutations of variables can be applied efficiently:
\meeq{
{\bfQ}_p(\rhos_n(g) \bfx)  = Q^\top {\bfP}_p(\rhos_n(g) \bfx)  = Q^\top \rho_{n,p}(g) {\bfP}_p(\bfx) = \sopmatrix{ \rho_1(g)^{\oplus a_1} \\ & \ddots \\ && \rho_r(g)^{\oplus a_r}} Q^\top {\bfP}_p(\bfx) \ccr
	= \sopmatrix{ \rho_1(g)^{\oplus a_1} \\ & \ddots \\ && \rho_r(g)^{\oplus a_r}} {\bfQ}_p(\bfx) 
}
 In \figref{polyblockdiag_p} we consider the practical implementation of this algorithm. Whilst the proven complexity is $\Ord(n d^4)$, in practice we observe  $\Ord(n d^3)$ complexity, resulting from dropping zero rows in the nullspace system. In \figref{cond_p} we measure the errors introduced by floating point arithmetic by comparing $Q^\top \rho_{n,p}(\tau_k) Q$ with the numerically computed $Q$ with the expected formula
 $$
 \sopmatrix{ \rho_1(g)^{\oplus a_1} \\ & \ddots \\ && \rho_r(g)^{\oplus a_r}}.
 $$
The error in this computation grows like $C d \epsilon_{\rm m}$ for machine epsilon $\epsilon_{\rm m} = 2^{-52} \approx 2.22 \times 10^{-16}$ where $C$ is a small constant, an error grow rate that is indicative of a numerical stable algorithm. 

}

\Figuretwo[cond_p,cond_n]
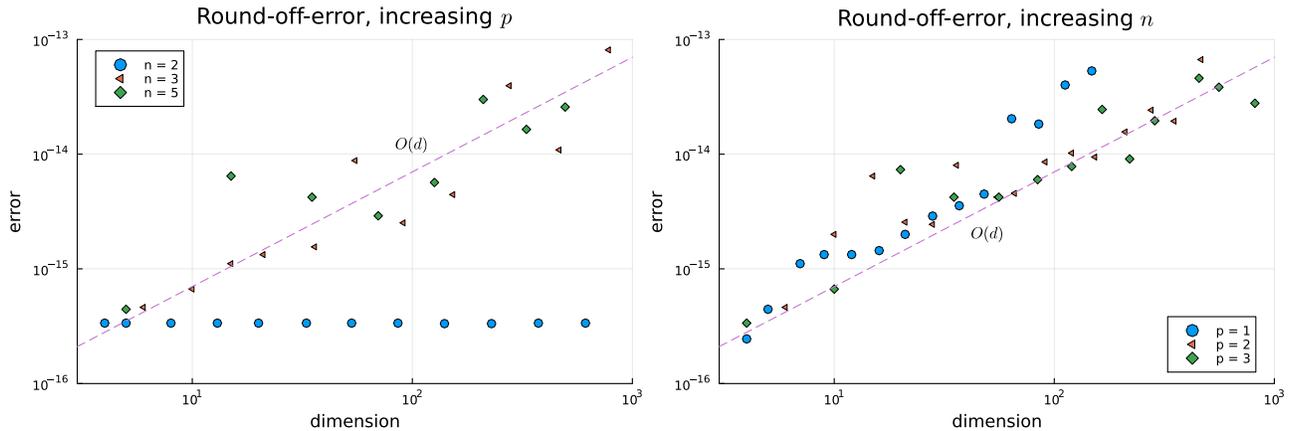
	{\addition The maximum entry-wise error comparing $Q^\top \rho_{n,p}(\tau_k) Q$ with the formula for the generators of the irreducible representations, where $Q$ is computed using the \algref{blockdiag} with floating point arithmetic. The error is proportional to the dimension of the representation which is indicative of a stable numerical algorithm. The left/right figures correspond to the same simulations as \figref{polyblockdiag_p}}.

{\addition 
\Section{future} Future work.

In this work we have established algorithms  for decomposing representations of $S_n$ that have polynomial complexity as the dimension of the representation or $n$ increase. These bounds are unlikely to be sharp, and in experiments  we often achieve  better complexity, though proving these sharper complexity results remains open and likely depends on a careful understanding of the sparsity present in the irreducible representations. The complexity can potentially be further reduced using randomised linear algebra. Whether the complexity can be reduced to such extent that general Kronecker coefficients can themselves be computed with high probability in polynomial time is unclear, and, given the super-algebraic growth in the dimensions of the representations, seems highly unlikely.

The practical implementation of the algorithm using floating point arithmetic is not guaranteed to be rigorous: it is possible round-off errors can alter the multiplicities or irreducible representations present. Fortunately, one can overcome this issue by employing techniques from the field of Validated Numerics \cite{tucker2011validated}.  This could potentially take the form of an implementation of the algorithm where the errors can be controlled rigorously by incorporating interval/ball arithmetic, see eg. \cite{Johansson2013arb} for an effective implementation of ball arithmetic that includes computing eigenvalues with rigorous bounds. More practically, the results of the current algorithm implemented with floating point arithmetic can be verified a postiori. To briefly sketch the argument: if $Q_{\rm fp} \approx Q$ is an approximation to a matrix that simultaneously diagonalises $\rho(X_j)$ (coming from \algref{multiplicities}, \algref{blockdiag}, or elsewhere) then we can write
$$
Q_{\rm fp} Q_{\rm fp}^\top = I + E
$$
where the norm of the matrix $E$ can be bounded using interval arithmetic.  Thus we can write
$$
Q_{\rm fp}^{-1} = Q_{\rm fp}^\top ( I + E)^{-1} = Q_{\rm fp}^\top + \tilde E
$$
where $\tilde E$ can be bounded explicitly using Neumann series in terms of $\| E \|$.  Thus we have
$$
Q_{\rm fp} \rho(X_j) Q_{\rm fp}^{-1} =  Q_{\rm fp} \rho(X_j) Q_{\rm fp}^\top +  \rho(X_j) \tilde E
$$
This matrix shares the same spectrum as $\rho(X_j)$ and them errors in the computation above can be bounded rigorously using interval arithmetic. The result will be a  matrix with small intervals surrounding integers on the diagonal combined with small intervals around 0 for off-diagonal entries. Gershgorin's theorem can then be combined with the fact that the eigenvalues must be integers to prove precisely which eigenvalues are present. To establish that the eigenvalues of each $\rho(X_j)$ share a common eigenvector and hence prove the correctness of $\Spec_\rho$ we can use the results of \cite{yamamoto1980error}.

The algorithms introduced may be useful for introducing sparsity in numerical discretisations. Representation theory has been used in the context of computing cubature rules for constructing sparse Hankel matrices associated with cubature rules on triangles, squares and hexagons   \cite{collowald2015moment}.  Other potential applications include solving  partial differential equations that commute with symmetry actions. Discretising such equations with  a basis associated with irreducible representations (such as the orthogonal polynomials constructed in \secref{examples})  will induce sparsity in the resulting discretisation. There is therefore a clear need for extension of the results of this paper for representations of general Coxeter groups, which are connected to orthogonal polynomials whose weight is invariant under symmetry actions (cf.~\cite{dunkl2014orthogonal}). Explicit irreducible representations are known for  the dihedral  \cite[Section 5.3]{serre1977linear} and the hyperoctohedral  \cite{geissinger1978representations,musili1993representations} groups which could potentially be used to adapt the proposed algorithm to these other groups.

}

\References

\ends